\documentstyle[12pt]{article}
\begin{document}
\title{WHY A CONJECTURE OF POINCARE DOESN'T WORK}
\author{B.G. Sidharth$^*$\\
Centre for Applicable Mathematics \& Computer Sciences\\
B.M. Birla Science Centre, Adarsh Nagar, Hyderabad - 500 063 (India)}
\date{}
\maketitle
\footnotetext{$^*$Email:birlasc@hd1.vsnl.net.in; birlard@ap.nic.in}
\begin{abstract}
Poincare had conjectured that the fact that closed loops could be shrunk
to points on a surface topologically equivalent to the surface of a
sphere can be generalised to three (and more) dimensions. After nearly a
century the conjecture has remained unproven. We given arguments below
to show that the conjecture doesn't work in three dimensions.
\end{abstract}
Nearly a hundred years ago, Poincare had conjectured that the fact that
closed loops could be shrunk to points on a two dimensional surface topologically
equivalent to the surface of a sphere can be generalised to three dimensions
also\cite{r1}. After all these years the conjecture has remained unproven.
We will now see why the three dimensional generalisation is not possible.\\
We firstly observe that a two dimensional surface on which closed smooth
loops can be shrunk continuously to arbitrarily small sizes is simply connected.
On such a surface we can define complex coordinates following the hydrodynamical
route exploiting the well known connection between the two. If we consider
laminar motion of an incompressible fluid we will have\cite{r2}
\begin{equation}
\vec \nabla \cdot \vec V = 0\label{e1}
\end{equation}
Equation (\ref{e1}) defines, as is well known, the stream function $\psi$ such
that
\begin{equation}
\vec V = \vec \nabla \times \psi \vec e_z\label{e2}
\end{equation}
where $\vec e_z$ is the unit vector in the $z$ direction.\\
Further, as the flow is irrotational, as well, we have
\begin{equation}
\vec \nabla \times \vec V = 0\label{e3}
\end{equation}
Equation (\ref{e3}) implies that there is a velocity potential $\phi$ such that,
\begin{equation}
\vec V = \vec \nabla \phi\label{e4}
\end{equation}
The equations (\ref{e2}) and (\ref{e4}) show that the functions $\psi$ and $\phi$
satisfy the Cauchy-Reimann equations of complex analysis\cite{r3}.\\
So it is possible to characterise the fluid elements by a complex variable
\begin{equation}
z = x + \imath y\label{e5}
\end{equation}
The question is can we generalise equation (\ref{e5}) to three dimensions? Infact
as Hamilton concluded a long time ago\cite{r4}, a generalisation leads not to
three but to four dimensions, with the three Pauli spin matrices $\vec \sigma$
replacing $\imath$. Further these Pauli spin matrices do not commute, and
characterise spin or vorticity.\\
This is not surprising - the reason lies in equation (\ref{e2}) or equivalently
in the multiplication law of complex numbers. (Infact, there is a general
tendency to loverlook this fact and this leads to the mistaken impression that
complex numbers are just an ordered pair of numbers, which latter are
usually associated with vectors.)\\
The above considerations give an explanation for the $3 + 1$ dimensionality
of space time\cite{r5}. Infact it is well known that it is the entanglement of
the spin of electrons (spin half) that implies three dimensionality of
physical space\cite{r6,r7}.

\end{document}